\documentclass[11pt]{article}
\usepackage{amsfonts}
\usepackage{amssymb}
\usepackage{amsmath}

\newcommand{\ignore}[1]{}
\ignore{The text in between these parentheses is ignored by LaTeX ignore}

\def\@begintheorem#1#2{\par\bgroup{\sc #1\ #2. }\it\ignorespaces}
\def\@opargbegintheorem#1#2#3{\par\bgroup{\sc #1\ #2\ (#3). } \it\ignorespaces}
\def\@endtheorem{\egroup}
\newtheorem{theorem}{Theorem}[section]
\newtheorem{lemma}[theorem]{Lemma}
\newtheorem{proposition}[theorem]{Proposition}
\newtheorem{example}[theorem]{Example}
\newcommand{\bt}[1]{\begin{theorem}\label{#1}}
\newcommand{\bl}[1]{\begin{lemma}\label{#1}}
\newcommand{\bp}[1]{\begin{proposition}\label{#1}}
\newcommand{\be}[1]{\begin{example}\rm\label{#1}}
\newcommand{\et}{\end{theorem}}
\newcommand{\el}{\end{lemma}}
\newcommand{\ep}{\end{proposition}}
\newcommand{\ee}{\end{example}}

\def \Z {\mathbb{Z}}
\def \C {{{\cal C}}}
\def \G {{{\cal G}}}
\def \type {{\rm type}}

\begin{document}

\title{\bf The Graver Complexity of Integer Programming}
\author{Yael Berstein \and Shmuel Onn}

\date{}
\maketitle

\begin{abstract}
In this article we establish an exponential lower bound on the
Graver complexity of integer programs. This provides new type of
evidence supporting the presumable intractability of integer programming.
Specifically, we show that the Graver complexity of the incidence
matrix of the complete bipartite graph $K_{3,m}$ satisfies
$g(m)=\Omega(2^m)$, with $g(m)\geq 17\cdot 2^{m-3}-7$ for every $m>3$ .
\vskip.2cm
\noindent
{\em keywords:} Graver basis, Gr\"obner basis, Graver complexity, Markov complexity,
contingency table, transportation polytope, transportation problem,
integer programming, computational complexity.

\end{abstract}

\section{Introduction}

In this article we establish an exponential lower bound on the
Graver complexity of integer programs. This provides new type of
evidence supporting the presumable intractability of integer programming.

We start by overviewing relevant recent developments in the theory
of integer programming which motivate our work and by providing
several definitions which are necessary for stating our result.

The {\em integer programming problem}, well known to be NP-complete, is to decide,
given integer $p\times q$ matrix $B$ and integer $p$ vector $b$ ,
if the following set of integer points in a polyhedron is nonempty,
$$S(B,b)\ :=\ \{x\in\Z^q\,:\,Bx=b\,,\ x\geq 0\}\ .$$
The {\em $n$-fold product} of an $s\times t$ matrix $A$ is the following
$(t+ns)\times nt$ matrix, with $I_t$ the $t\times t$ identity:

$$A^{(n)}\quad:=\quad ({\bf 1}_n\otimes I_t)\oplus (I_n \otimes A)\quad=\quad
\left(
\begin{array}{ccccc}
  I_t    & I_t    & I_t    & \cdots & I_t    \\
  A  & 0      & 0      & \cdots & 0      \\
  0  & A      & 0      & \cdots & 0      \\
  \vdots & \vdots & \ddots & \vdots & \vdots \\
  0  & 0      & 0      & \cdots & A      \\
\end{array}
\right)\quad .
$$
The efficient solution of {\em $n$-fold integer programming}
in variable dimension was recently proved in \cite{DHOW}:
\bp{N-fold}
Fix any integer matrix $A$. Then there is a polynomial time algorithm
that, given $n$ and integer $(t+ns)$ vector $a$, decides if the set
$S(A^{(n)},a)=\{x\in\Z^{nt}:A^{(n)}x=a,\ x\geq 0\}$ is nonempty.
\ep
The time complexity of the algorithm underlying Proposition \ref{N-fold} is
$O\left(n^{g(A)}\|\log|a|\|_1\right)$ where $g(A)$ is the {\em Graver complexity}
of the matrix $A$. We proceed to define this notion, recently introduced in \cite{SS}.

Define a partial order $\sqsubseteq$ on $\Z^n$ which extends the coordinate-wise
order $\leq$ on $\Z_+^n$ as follows: for two vectors $u,v\in\Z^n$ put
$u\sqsubseteq v$ if $|u_i|\leq |v_i|$ and $u_iv_i\geq 0$ for $i=1,\ldots,n$.
A suitable extension of the classical lemma of Gordan \cite{Gor} implies
that every subset of $\Z^n$ has finitely-many $\sqsubseteq$-minimal elements.
The {\em Graver basis} of an integer matrix $A$, introduced in \cite{Gra}, is
defined to be the finite set $\G(A)$ of $\sqsubseteq$-minimal elements in
the set $\{x\in\Z^n:Ax=0\,,\ x\neq 0\}$ of nontrivial integer dependencies on $A$.

Consider any $s\times t$ integer matrix $A$. For any positive integer $n$ consider
vectors $x\in\Z^{nt}$ indexed as $x=(x^1,\dots,x^n)$ with each block $x^i$
lying in $\Z^t$. The {\em type} of $x=(x^1,\dots,x^n)$ is the number
$\type(x):=|\{i\,:\,x^i\neq 0\}|$ of nonzero blocks of $x$.
The {\em Graver complexity} of $A$ is defined to be
$$g(A)\ :=\ \sup\left(\{0\}\cup
\left\{\type(x)\ :\ x\in \bigcup_{n\geq 1}\G(A^{(n)})\right\}\right)\ .$$
(If the columns of $A$ are linearly independent then
$\G(A^{(n)})$ is empty for all $n$ and hence $g(A)=0$.)
The following result was recently proved in \cite{SS}, extending a result
of \cite{AT} for the matrix in (\ref{matrix}) below:
\bp{GraverComplexity}
The Graver complexity $g(A)$ of every integer matrix $A$ is finite.
\ep

Let $(1,1,1)^{(m)}$ be the $m$-fold product of the $1\times 3$ matrix
$(1,1,1)$. Note that $(1,1,1)^{(m)}$ is precisely the $(3+m)\times 3m$ vertex-edge
incidence matrix of the complete bipartite graph $K_{3,m}$. For instance,
\begin{equation}\label{matrix}
(1,1,1)^{(3)}\ =\
\left(
\begin{array}{ccccccccc}
  1 & 0 & 0 & 1 & 0 & 0 & 1 & 0 & 0 \\
  0 & 1 & 0 & 0 & 1 & 0 & 0 & 1 & 0 \\
  0 & 0 & 1 & 0 & 0 & 1 & 0 & 0 & 1 \\
  1 & 1 & 1 & 0 & 0 & 0 & 0 & 0 & 0 \\
  0 & 0 & 0 & 1 & 1 & 1 & 0 & 0 & 0 \\
  0 & 0 & 0 & 0 & 0 & 0 & 1 & 1 & 1 \\
\end{array}
\right)\ .
\end{equation}
A recent universality theorem in \cite{DO1} asserts that {\em every} bounded set $S(B,b)$
stands in polynomial-time computable linear bijection with the set of integer points
$S\left(\left((1,1,1)^{(m)}\right)^{(n)},a\right)$ for some $m$, $n$ and $a$:
\bp{Universality}
There is a polynomial time algorithm that, given $B$ and $b$
with $S(B,b)$ bounded, computes $m$, $n$, and integer
$(3m+n(3+m))$ vector $a$ such that $S(B,b)$ stands in linear bijection with
$$S\left(\left((1,1,1)^{(m)}\right)^{(n)},a\right)\ =\
\left\{x\in\Z^{3mn}\ :\ \left((1,1,1)^{(m)}\right)^{(n)} x\, =\, a\,,\ \ x\geq 0\right\}\ .$$
\ep

Let $g(m):=g\left((1,1,1)^{(m)}\right)$ denote the Graver complexity of $(1,1,1)^{(m)}$.
Proposition \ref{N-fold} and Proposition \ref{Universality} then imply the following
interestingly contrasting situations about the computational complexity of deciding
if $S\left(\left((1,1,1)^{(m)}\right)^{(n)},a\right)$ is nonempty:
for every fixed $m$, the problem is decidable in polynomial time
$O\left(n^{g(m)}\|\log|a|\|_1\right)$; but for variable $m$, the problem is NP-complete.
So if $P\neq NP$ then $g(m)$ cannot be bounded by a constant and
must grow as a function of $m$. In this article we show that, in fact, it grows
{\em exponentially fast}, as $g(m)=\Omega(2^m)$. We establish the following statement.
\bt{Exponential}
The Graver complexity of $(1,1,1)^{(m)}$ satisfies
$g(m)\geq 17\cdot 2^{m-3}-7$ for every $m>3$.
\et

Theorem \ref{Exponential} implies the exponential lower bound $g(m)=\Omega(2^m)$.
An exponential upper bound $g(m)=O\left(m^4 6^m\right)$ can be derived from
the Cramer rule and the Hadamard bound. Narrowing the gap between these
bounds remains a challenging and important problem.
One possible approach might be to make a careful use of the complex
universality constructions of \cite{DO1} and \cite{DO2}.
Note that computing $g(m)$ even for small $m$ is extremely difficult:
while it is known that $g(3)=9$,
the exact value of $g(m)$ is unknown for all $m>3$.
The lower bound provided by Theorem \ref{Exponential} is the sharpest one known for every $m>3$.
We point out that Theorem \ref{Exponential} implies an exponential lower bound also
on the {\em Gr\"obner complexity} of $(1,1,1)^{(m)}$ considered recently in \cite{HN}
as well as on the {\em Markov complexity} of $(1,1,1)^{(m)}$.

We conclude by introducing a new graph invariant that naturally arises
in this context and deserves further study. Let $G=(V,E)$ be a graph or a digraph and
let $A$ be its $V\times E$ incidence matrix. (For a graph, $A_{v,e}$ is $1$ if vertex $v$
lies in edge $e$ and is $0$ otherwise; for a digraph, $A_{v,e}$ is $1$ if vertex $v$
is the head of arc $e$, is $-1$ if $v$ is the tail of $e$, and is $0$ otherwise.)
Define the {\em Graver complexity} of $G$ to be the Graver complexity of its incidence
matrix $A$, that is, $g(G):=g(A)$. In particular, since $(1,1,1)^{(m)}$ is the
incidence matrix of the complete bipartite graph $K_{3,m}$, it follows that
$g(K_{3,m})=g(m)$ is precisely the function studied here.
It is quite intriguing that the Graver complexity of $K_{3,4}$ is yet unknown.

Returning to integer programming, the $n$-fold integer programming problem associated
with the incidence matrix of a graph or a digraph $G$ is the corresponding
{\em $n$-commodity $b$-matching problem} or {\em $n$-commodity transshipment problem}
over $G$, respectively. The Graver complexity $g(G)$ controls the computational complexity
of solving these problems over $G$. These problems will be studied elsewhere.

\section{Proof}

Our starting point is the following characterization of the Graver
complexity from \cite{SS}. Here $\G(\G(A))$ denotes the Graver basis of a
matrix whose columns are the elements of $\G(A)$ ordered arbitrarily.
\bp{Characterization}
The Graver complexity of every $A$ satisfies
$g(A)=\max\left\{\|x\|_1:x\in \G(\G(A))\,\right\}$.
\ep

A {\em circuit} of an integer matrix $A$ is a nonzero integer vector $x$ satisfying $Ax=0$,
that has inclusion-minimal support with respect to this property, and whose nonzero entries
are relatively prime. Let $\C(A)$ denote the (finite) set of circuits of the matrix $A$.
The following statement is well known, see \cite{Stu}.

\bp{Circuits}
The set of circuits and the Graver basis of every $A$ satisfy $\C(A)\subseteq \G(A)$.
\ep

A linear relation $\sum_{i=1}^k h_i v^i=0$ on integer vectors $v^1,\dots,v^k$
is {\em primitive} if the coefficients $h_1,\dots,h_k$ are relatively prime
positive integers and no $k-1$ of the $v^i$ satisfy any nontrivial linear relation.
Our interest in circuits and primitive relations stems from the following statement.

\bp{Circuit}
Suppose that $\sum_i h_i x^i=0$ is a primitive relation on some circuits $x^i$ of
$(1,1,1)^{(m)}$. Then the Graver complexity of $(1,1,1)^{(m)}$ satisfies $g(m)\geq \sum_i h_i$.
\ep
{\bf \em Proof.}$\quad$
By Proposition \ref{Circuits} we have the containment
$\C\left((1,1,1)^{(m)}\right) \subseteq \G\left((1,1,1)^{(m)}\right)$,
and therefore the circuits $x^i$ appear among the elements of the Graver basis
of $(1,1,1)^{(m)}$. Let $G$ be a matrix whose columns are the elements of
$\G\left((1,1,1)^{(m)}\right)$ with the $x^i$ coming first.
Then the vector $h$ of dimension $\left|\G\left((1,1,1)^{(m)}\right)\right|$,
that consists of the coefficients $h_i$ of the given relation augmented with
sufficiently many trailing zeros, is a circuit of $G$.
By Proposition \ref{Circuits} applied once more, we find that
$$h\ \in\ \C(G)\ \subseteq\ \G(G)\ =\ \G\left(\G\left((1,1,1)^{(m)}\right)\right)\ .$$
The claim then follows since by Proposition \ref{Characterization} we have
$$g(m)\ =\ \max\left\{\|x\|_1\ :\ x\in \G\left(\G\left((1,1,1)^{(m)}\right)\right)\,\right\}
\ \geq\ \|h\|_1\ =\ \sum_i h_i\ . \quad\quad\blacksquare$$

We employ below the following notation, where $m$ is
any integer understood from the context. Let
$$A\,:=\,\{a,b,c\},\ \ U\,:=\,\{u_1,\dots,u_m\},
\ \ V\,:=\,A\uplus U,\ \ E\,:=\,A\times U\ .$$
Then $V$ and $E$ are, respectively, the set of vertices and set of edges
of the complete bipartite graph $K_{3,m}$, and index, respectively, the
rows and columns of its vertex-edge incidence matrix $(1,1,1)^{(m)}$.
It will be convenient to interpret each vector $x\in\Z^E$ also as:
(1) an integer valued function on the set of edges $E=A\times U$;
(2) a $3\times m$ matrix with rows and columns indexed, respectively, by $A$ and $U$.
With these interpretations, $x$ is in $\C\left((1,1,1)^{(m)}\right)$ if and only if:
(1) as a function on $E$, its support is a circuit of $K_{3,m}$, along which it
alternates in values $\pm 1$, and can be indicated by the sequence $(v_1,v_2,\dots,v_l)$
of vertices of the circuit of $K_{3,m}$ on which it is supported, with the
convention that its value is $+1$ on the first edge $(v_1,v_2)$ in that sequence;
(2) as a matrix, it is nonzero, has $0,\pm 1$ entries, has zero row and column sums,
and has inclusion-minimal support with respect to these properties.

Here is an example, that will also play a role in the proofs below,
demonstrating this notation.

\be{Example}{\bf (A lower bound for m\ =\ 4).}\ \
Let $m=4$. Define the following seven circuits of $K_{3,4}$:
$$
\begin{tabular}{ll|cccc|c}
  && $u_1$ & $u_2$ & $u_3$ & $u_4$ & \\
 \hline
  &&  $0$ & $0$ & $\!\!\!\!-1$ & $1$ & $a$  \\
$x^1\ :=\ (a,u_4,c,u_2,b,u_3)$ & = &
      $0$ & $\!\!\!\!-1$ &   $1$   & $0$  & $b$  \\
  &&  $0$ & $1$ &   $0$   & $\!\!\!\!-1$  & $c$  \\
 \hline
  &&  $\!\!\!\!-1$ & $1$ & $0$ & $0$ & $a$  \\
$x^2\ :=\ (a,u_2,c,u_3,b,u_1)$ & = &
      $1$ & $0$ &   $\!\!\!\!-1$   & $0$  & $b$  \\
  &&  $0$ & $\!\!\!\!-1$ &   $1$   & $0$  & $c$  \\
 \hline
  &&  $0$ & $\!\!\!\!-1$ & $0$ & $1$ & $a$  \\
$x^3\ :=\ (a,u_4,b,u_1,c,u_2)$ & = &
      $1$ & $0$ &   $0$   & $\!\!\!\!-1$  & $b$  \\
  &&  $\!\!\!\!-1$ & $1$ &   $0$   & $0$  & $c$  \\
 \hline
  &&  $\!\!\!\!-1$ & $0$ & $0$ & $1$ & $a$  \\
$x^4\ :=\ (a,u_4,b,u_2,c,u_1)$ & = &
      $0$ & $1$ &   $0$   & $\!\!\!\!-1$  & $b$  \\
  &&  $1$ & $\!\!\!\!-1$ &   $0$   & $0$  & $c$  \\
 \hline
  &&  $1$ & $0$ & $\!\!\!\!-1$ & $0$ & $a$  \\
$x^5\ :=\ (a,u_1,b,u_2,c,u_3)$ & = &
      $\!\!\!\!-1$ & $1$ &   $0$   & $0$  & $b$  \\
  &&  $0$ & $\!\!\!\!-1$ &   $1$   & $0$  & $c$  \\
 \hline
  &&  $0$ & $\!\!\!\!-1$ & $1$ & $0$ & $a$  \\
$x^6\ :=\ (a,u_3,b,u_4,c,u_2)$ & = &
      $0$ & $0$ &   $\!\!\!\!-1$   & $1$  & $b$  \\
  &&  $0$ & $1$ &   $0$   & $\!\!\!\!-1$  & $c$  \\
 \hline
  &&  $0$ & $1$ & $0$ & $\!\!\!\!-1$ & $a$  \\
$x^7\ :=\ (a,u_2,b,u_3,c,u_4)$ & = &
      $0$ & $\!\!\!\!-1$ &   $1$   & $0$  & $b$  \\
  &&  $0$ & $0$ &   $\!\!\!\!-1$   & $1$  & $c$  \\
 \hline
\end{tabular}
$$

\noindent
Then the circuits $x^i$ satisfy the primitive relation $x^1+2x^2+3x^3+3x^4+5x^5+6x^6+7x^7=0$.
Therefore, by Proposition \ref{Circuit} we obtain the bound
$g(4)=g\left((1,1,1)^{(4)}\right)\geq 1+2+3+3+5+6+7=27$.
\ee
We have the following lemma.
\bl{Main}
Suppose there are $k$ circuits $x^i$ of $(1,1,1)^{(m)}$ admitting a primitive
relation $\sum_i h_i x^i=0$ with $x^k=(a,u_{m-2},b,u_{m-1},c,u_m)$ and $h_k$ odd.
Then there are $k+2$ circuits ${\bar x}^i$ of $(1,1,1)^{(m+1)}$ admitting
a primitive relation $\sum_i{\bar h}_i {\bar x}^i=0$ with
${\bar x}^{k+2}=(a,u_{m-1},b,u_m,c,u_{m+1})$ and ${\bar h}_{k+2}$ odd, where
\begin{equation}\label{Induction}
{\bar h}_i\ =\ 2h_i\,,\quad i=1,\dots,k-1\,,\quad\quad
{\bar h}_{k+2}\ =\ {\bar h}_{k+1}\ =\ {\bar h}_k\ =\ h_k\ .
\end{equation}
\el
{\bf \em Proof.}$\quad$
Using the natural embedding of the complete bipartite graph $K_{3,m}$ into $K_{3,m+1}$,
we can interpret circuits of the former also as circuits of the latter.
Put $y^i:=x^i$ for $i=1,\dots,k-1$ and define
$$
\begin{tabular}{ll|cccccc|c}
  &&  $u_1$ & $\cdots$ & $u_{m-2}$ & $u_{m-1}$ & $u_m$ & $u_{m+1}$ & \\
 \hline
  &&  $0$ & $\cdots$ & $1$ & $0$ & $0$ & $\!\!\!\!-1$ & $a$  \\
$y^k\ :=\ (a,u_{m-2},b,u_{m-1},c,u_{m+1})$ & = &
      $0$ & $\cdots$ & $\!\!\!\!-1$ & $1$ &   $0$   & $0$  & $b$  \\
  &&  $0$ & $\cdots$ & $0$ & $\!\!\!\!-1$ &   $0$   & $1$  & $c$  \\
 \hline
  &&  $0$ & $\cdots$ & $1$ & $0$ & $\!\!\!\!-1$ & $0$ & $a$  \\
$y^{k+1}\ :=\ (a,u_{m-2},b,u_{m+1},c,u_m)$ & = &
      $0$ & $\cdots$ & $\!\!\!\!-1$ & $0$ &   $0$   & $1$  & $b$  \\
  &&  $0$ & $\cdots$ & $0$ & $0$ &   $1$   & $\!\!\!\!-1$  & $c$  \\
 \hline
  &&  $0$ & $\cdots$ & $0$ & $0$ & $\!\!\!\!-1$ & $1$ & $a$  \\
$y^{k+2}\ :=\ (a,u_{m+1},b,u_{m-1},c,u_m)$ & = &
      $0$ & $\cdots$ & $0$ & $1$ &   $0$   & $\!\!\!\!-1$  & $b$  \\
  &&  $0$ & $\cdots$ & $0$ & $\!\!\!\!-1$ &   $1$   & $0$  & $c$  \\
 \hline
\end{tabular}
$$
Note that these circuits satisfy $y^k+y^{k+1}+y^{k+2}=2x^k$.
Suppose that $\sum_i{\bar h}_i y^i=0$ is a nontrivial relation on the $y^i$.
Without loss of generality we may assume that the ${\bar h}_i$ are relatively prime
integers, at least one of which is positive.
Since the edges $(a,u_{m+1}),(b,u_{m+1}),(c,u_{m+1})$ of $K_{3,m+1}$ are
not in $K_{3,m}$ and hence in no circuit $y^i$ for $i<k$, the restrictions of the
relation $\sum_i{\bar h}_i y^i=0$ to these edges (or to the corresponding matrix entries)
forces the equalities ${\bar h}_{k+2}={\bar h}_{k+1}={\bar h}_k$. We then obtain
$$0\ =\ \sum_{i=1}^{k+2}{\bar h}_i y^i\ =\ \sum_{i=1}^{k-1}{\bar h}_i y^i
\ +\ {\bar h}_k (y^k+y^{k+1}+y^{k+2})\ =\ \sum_{i=1}^{k-1}{\bar h}_i x^i
\ +\ {\bar h}_k (2x^k)\ =\ \sum_{i=1}^{k-1}{\bar h}_i x^i \ +\ 2{\bar h}_k x^k$$
which is a nontrivial integer relation on the $x^i$. So there must exist an integer
$\alpha$ so that, for all $i$, the coefficient of $x^i$ in that relation is
$\alpha$ times the coefficient of $x^i$ in the relation $\sum_{i=1}^k h_i x^i=0$, that is,
$${\bar h}_i\ =\ \alpha h_i\,,\quad i=1,\dots,k-1\,,\quad\quad 2{\bar h}_k\ =\ \alpha h_k\ .$$
Since all the $h_i$ and at least one of the ${\bar h}_i$ are positive, these equations imply
that $\alpha$ is positive. Therefore all ${\bar h}_i$ are positive, implying that
the relation $\sum_i{\bar h}_i y^i=0$ on the $y^i$ is primitive.
Since $h_k$ is odd, the equation $2{\bar h}_k=\alpha h_k$ implies
that $\alpha$ is even and therefore $\alpha=2\mu$ for some positive integer $\mu$,
implying ${\bar h}_k=\mu h_k$. Then $\mu$ divides each of the ${\bar h}_i$, which are relatively
prime, and therefore $\mu=1$ and $\alpha=2$. It follows that the ${\bar h}_i$ satisfy
equation (\ref{Induction}) and in particular ${\bar h}_{k+2}=h_k$ is odd as claimed.

Now apply to the vertices of $K_{3,m+1}$ a permutation that maps $u_{m+1},u_{m-1},u_m$
to $u_{m-1},u_m,u_{m+1}$ in that order and fixes the rest of the vertices.
For $i=1,\dots,k+2$ let ${\bar x}^i$ be the circuit of $K_{3,m+1}$ that is
the image of $y^i$ under this permutation. Then the ${\bar x}^i$ also satisfy the
primitive relation $\sum_i{\bar h}_i {\bar x}^i=0$ with the same coefficients ${\bar h}_i$,
and ${\bar x}^{k+2}=(a,u_{m-1},b,u_m,c,u_{m+1})$. This completes the proof.
$\quad\quad\blacksquare$

We are now in position to prove our theorem.

\vskip0.5cm

{\bf \em Proof of Theorem \ref{Exponential}.}$\quad$
We prove by induction on $m$ that, for all $m\geq 4$, there are $2m-1$
circuits $x^i$ of $(1,1,1)^{(m)}$ with $x^{2m-1}=(a,u_{m-2},b,u_{m-1},c,u_m)$,
satisfying a primitive relation $\sum_i h_i x^i=0$ with $h_{2m-1}=7$ and
$\sum_i h_i=17\cdot 2^{m-3}-7$.
This combined with Proposition \ref{Circuit} implies the theorem.

The basis of the induction, at $m=4$, is verified by the seven
circuits constructed in Example \ref{Example}.

Suppose now the hypothesis holds for some $m\geq 4$ and let $x^i$ be $2m-1$ circuits
with corresponding coefficients $h^i$ verifying the hypothesis. Lemma \ref{Main}
applied to this data with $k=2m-1$ then guarantees the existence of
$k+2=2m+1=2(m+1)-1$ circuits ${\bar x}^i$ with corresponding coefficients
${\bar h}^i$ satisfying $x^{2(m+1)-1}=(a,u_{m-1},b,u_m,c,u_{m+1})$ and
${\bar h}_{2(m+1)-1}={\bar h}_{2m}={\bar h}_{2m-1}=h_{2m-1}=7$, and, moreover,
$$\sum_{i=1}^{2m+1}{\bar h}_i\ =\ \sum_{i=1}^{2m-2} 2 h_i\, +\, 3 h_{2m-1}
\ =\ 2\sum_{i=1}^{2m-1} h_i\, +\, h_{2m-1}\ =\ 2\left(17\cdot 2^{m-3}-7\right)
\, +\,  7 \ =\ 17\cdot 2^{(m+1)-3}-7\ .\quad\ \ \blacksquare$$

\section{Example}

We conclude by exhibiting nine circuits of $(1,1,1)^{(5)}$ and $K_{3,5}$,
obtained by applying our construction,
$$x^1\ =\ (a,u_5,c,u_2,b,u_4)\,,\quad  x^2\ =\ (a,u_2,c,u_4,b,u_1)\,,\quad
  x^3\ =\ (a,u_5,b,u_1,c,u_2)\ ,$$
$$x^4\ =\ (a,u_5,b,u_2,c,u_1)\,,\quad  x^5\ =\ (a,u_1,b,u_2,c,u_4)\,,\quad
  x^6\ =\ (a,u_4,b,u_5,c,u_2)\ ,$$
$$x^7\ =\ (a,u_2,b,u_4,c,u_3)\,,\quad  x^8\ =\ (a,u_2,b,u_3,c,u_5)\,,\quad
  x^9\ =\ (a,u_3,b,u_4,c,u_5)\ ,$$
that satisfy the primitive relation
$$2x^1\ +\ 4x^2\ +\ 6x^3\ +\ 6x^4\ +\ 10x^5
\ +\ 12x^6\ +\ 7x^7\ +\ 7x^8\ +\ 7x^9\quad =\quad 0\ ,$$
thereby demonstrating the lower bound $g(5)\geq 61$ on the
Graver complexity of $(1,1,1)^{(5)}$ and $K_{3,5}$.

\section*{Acknowledgements}

The research of Yael Berstein was partially supported by an Irwin and Joan
Jacobs Scholarship and by a scholarship from the Technion Graduate School.
The research of Shmuel Onn was partially supported by the ISF - Israel Science
Foundation and by the Fund for the Promotion of Research at the Technion.

\noindent {\small Yael Berstein}\newline
\emph{Technion - Israel Institute of Technology, 32000 Haifa, Israel}\newline
\emph{email: yaelber{\small @}tx.technion.ac.il}

\vskip.3cm\noindent {\small Shmuel Onn}\newline
\emph{Technion - Israel Institute of Technology, 32000 Haifa, Israel}\newline
\emph{email: onn{\small @}ie.technion.ac.il},
\ \ \emph{http://ie.technion.ac.il/{\small $\sim$onn}}

\end{document}